\documentclass{ifacconf}

\usepackage{graphicx}
\usepackage{natbib}

\usepackage{amsmath}
\usepackage{amssymb}
\usepackage{mathtools}
\usepackage{enumitem}
\usepackage{xcolor}
\usepackage{booktabs}
\usepackage{tikz}
\usetikzlibrary{arrows, arrows.meta}
\usepackage[european, nooldvoltagedirection]{circuitikz}

\DeclareMathOperator{\Real}{Re} 
\DeclareMathOperator{\Imag}{Im} 
\DeclareMathOperator{\id}{Id} 
\DeclareMathOperator{\diff}{d\!} 

\DeclarePairedDelimiter{\abs}{\lvert}{\rvert} 

\setlist[enumerate]{label={\textnormal{(\alph*)}}, ref={(\alph*)}, leftmargin=*, nosep}

\newcommand{\suchthat}{\ifnum\currentgrouptype=16 \mathrel{}\middle|\mathrel{}\else\mid\fi}
\newcommand{\transpose}{^{\mathrm T}}
\begin{document}
\begin{frontmatter}

\title{Effects of Roots of Maximal Multiplicity on the Stability of Some Classes of Delay Differential-Algebraic Systems: The Lossless Propagation Case\thanksref{footnoteinfo}}

\thanks[footnoteinfo]{Corresponding author: Guilherme Mazanti (gui\-lher\-me.ma\-zan\-ti@l2s.cen\-trale\-sup\-elec.fr). G.M., I.B., and S.I.N.\ are also with Inria, DISCO Team, Saclay -- \^Ile-de-France research center.}

\author[IPSA,LSS]{Guilherme Mazanti} 
\author[IPSA,LSS]{Islam Boussaada} 
\author[LSS]{Silviu-Iulian Niculescu}
\author[LSS]{Yacine Chitour}

\address[LSS]{Universit\'e Paris-Saclay, CNRS, CentraleSup\'elec, Laboratoire des signaux et syst\`emes, 91190, Gif-sur-Yvette, France.\\\{first name.last name\}@l2s.centralesupelec.fr}
\address[IPSA]{Institut Polytechnique des Sciences Avanc\'ees (IPSA)\\63 boulevard de Brandebourg, 94200 Ivry-sur-Seine, France.}

\begin{abstract}
It has been observed in several recent works that, for some classes of linear time-delay systems, spectral values of maximal multiplicity are dominant, a property known as multiplicity-induced-dominancy (MID). This paper starts the investigation of whether MID holds for delay differential-algebraic systems by considering a single-delay system composed of two scalar equations. After motivating this problem and recalling some recent results for retarded delay differential equations, we prove that the MID property holds for the delay differential-algebraic system of interest and present some applications.
\end{abstract}

\begin{keyword}
Time-delay equations, stability analysis, spectral methods, root assignment.
\end{keyword}

\end{frontmatter}

\section{Introduction}

Time delays are useful modeling tools in a wide range of scientific and technological domains such as biology, chemistry, economics, physics, or engineering. They may represent, for instance, the reaction time of an engineering system, the transfer time of material, energy, or information between parts of a system, the duration of a chemical reaction, or the duration of maturation processes in biology. Due to these applications and the challenging mathematical problems arising in their analysis, systems with time delays have been the subject of much attention by researchers in several fields, in particular since the 1950s and 1960s, such as, for instance, in \cite{Bellman1963Differential, Halanay1966Differential, Pinney1958Ordinary}. We refer to \cite{Diekmann1995Delay, Gu2003Stability, Hale1993Introduction, Insperger2011Semi, Michiels2007Stability, Stepan1989Retarded} for details on time-delay systems and their applications.

This paper is interested in the analysis of stability and asymptotic behavior of systems with delays of the general form
\begin{equation}
\label{IntroDiffAlgeb}
\left\{
\begin{aligned}
x^\prime(t) & = A x(t) + \sum_{k=1}^N B_k y(t - \tau_k), \\
y(t) & = C x(t) + \sum_{k=1}^N D_k y(t - \tau_k),
\end{aligned}
\right.
\end{equation}
where $x(t) \in \mathbb R^{d_x}$, $y(t) \in \mathbb R^{d_y}$, $N$, $d_x$, and $d_y$ are positive integers, $\tau_1, \dotsc, \tau_N$ are positive delays, and, for $k \in \{1, \dotsc, N\}$, $A$, $B_k$, $C$, and $D_k$ are matrices of appropriate dimensions. Systems such as \eqref{IntroDiffAlgeb} are \emph{delay differential-algebraic systems} since they are written as a system of delay differential equations coupled with a system of algebraic equations with delays. They correspond to a particular class of delay differential-algebraic systems in which delays appear only in the $y$ variable, which is the case in particular in models of lossless propagation phenomena (see Section~\ref{SecMotivation} for an example).

Differential-algebraic systems have been extensively studied in the delay-free setting (see, e.g., \cite{Brenan1989Numerical, Coleman1998Differential, Griepentrog1986Differential, Kumar1999Control, Kunkel2006Differential}). This kind of system arises naturally in several situations, such as in some electronic circuit models, in some control problems with constraints, or in the limiting behavior of singularly perturbed systems. Delay differential-algebraic systems have also been considered in the literature, arising in general from systems of transport partial differential equations representing some propagation phenomenon and coupled with static and dynamic boundary conditions (see, e.g., \cite{Niculescu2001Delay, Halanay1997Stability, Hale1993Introduction, Brayton1968Small}). More specific motivating examples are described in Section~\ref{SecMotivation}.

The stability analysis of time-delay systems has attracted much research effort and is an active field (see, e.g., \cite{Cooke1986Zeroes, Gu2003Stability, Michiels2007Stability, Olgac2002Exact, Sipahi2011Stability}). Similarly to the delay-free situation, one may address the asymptotic behavior of a linear time-invariant time-delay system through spectral methods by considering the corresponding characteristic function, whose complex roots determine the asymptotic behavior of solutions of the system (see, e.g., \cite{Hale1993Introduction, Michiels2007Stability, Mori1982Estimate}). The characteristic function of \eqref{IntroDiffAlgeb} is the function $\Delta: \mathbb C \to \mathbb C$ defined for $s \in \mathbb C$ by
\begin{equation}
\label{IntroDelta}
\Delta(s) = \det\left(s E - \widehat A - \sum_{k=1}^N e^{-s \tau_k} \widehat B_k\right),
\end{equation}
where $E, \widehat A, \widehat B_1, \dotsc, \widehat B_N$ are $(d_x + d_y) \times (d_x + d_y)$ matrices defined by blocks as
\[
E = \begin{pmatrix}\id_{d_x} & 0 \\ 0 & 0 \\\end{pmatrix}, \quad \widehat A = \begin{pmatrix}A & 0 \\ C & -\id_{d_y}\\\end{pmatrix}, \quad \widehat B_k = \begin{pmatrix}0 & B_k\\ 0 & D_k\\\end{pmatrix},
\]
for $k \in \{1, \dotsc, N\}$. Notice that \eqref{IntroDiffAlgeb} can be rewritten in terms of these matrices as
\[
E z^\prime(t) = \widehat A z(t) + \sum_{k=1}^N \widehat B_k z(t - \tau_k),
\]
where $z(t) = (x(t)\transpose, y(t)\transpose)\transpose$. Similarly to the delay-free case, the exponential behavior of \eqref{IntroDiffAlgeb} is determined by the \emph{spectral abscissa} $\gamma$ of $\Delta$, defined by $\gamma = \sup\{\Real s \mid s \in \mathbb C \text{ and } \Delta(s) = 0\}$, and all solutions of \eqref{IntroDiffAlgeb} converge exponentially to $0$ if and only if $\gamma < 0$.

The spectral abscissa of $\Delta$ is related to the notion of dominant root, defined as follows.
\begin{defn}
Let $Q: \mathbb C \to \mathbb C$ and $s_0 \in \mathbb C$ be such that $Q(s_0) = 0$. We say that $s_0$ is a \emph{dominant} (respectively, \emph{strictly dominant}) root of $Q$ if, for every $s \in \mathbb C \setminus \{s_0\}$ such that $Q(s) = 0$, one has $\Real s \leq \Real s_0$ (respectively, $\Real s < \Real s_0$).
\end{defn}
It follows immediately from the above definition that, if $Q$ admits a dominant root $s_0$, then $\gamma = \Real s_0$, but dominant roots may not exist in general.

Functions of the form \eqref{IntroDelta} are particular instances of quasipolynomials, whose definition is the following.
\begin{defn}
\label{DefiQuasipolynomial}
A \emph{quasipolynomial} is an entire function $Q$ which can be written under the form
\[
Q(s) = \sum_{k=0}^\ell P_k(s) e^{\lambda_k s},
\]
where $\ell$ is a positive integer, $\lambda_0, \dotsc, \lambda_\ell$ are pairwise distinct real numbers, and, for $k \in \{0, \dotsc, \ell\}$, $P_k$ is a non-zero polynomial of degree $d_k$. The integer $D = \ell + \sum_{k=0}^\ell d_k$ is called the \emph{degree} of $Q$.
\end{defn}
The above definition of the degree of a quasipolynomial is motivated by a classical property, provided in \cite[Problem 206.2]{Polya1998Problems} and known as the \emph{P\'{o}lya--Szeg\H{o} bound}, which implies that, given a quasipolynomial $Q$ of degree $D \geq 0$, the multiplicity of any root of $Q$ does not exceed $D$. Recent works such as \cite{Boussaada2016Characterizing, Boussaada2016Tracking} have provided characterizations of multiple roots of quasipolynomials using approaches based on Birkhoff and Vandermonde matrices.

It has been recently remarked (see, e.g., \cite{Boussaada2016Tracking, Boussaada2018Further, Boussaada2020Multiplicity, MazantiQualitative, MazantiSpectral}) that, for quasipolynomials coming from some systems with time-delays, real roots of maximal multiplicity are often dominant, a property usually referred to as \emph{multiplicity-induced-dominancy} (MID for short). This property has been shown to hold, in particular, for scalar single-delay differential equations of retarded type (see \cite{MazantiQualitative}), a result we present and explain below in Section~\ref{SecRetarded}, and also extended for some systems to the case of complex roots of maximal multiplicity in \cite{MazantiSpectral}.

One of the applications of the MID property is in the design of stabilizing feedback controllers for control systems with time delays, as in \cite{Boussaada2020Multiplicity}. A major difficulty when addressing this question is that, except in degenerate situations, quasipolynomials have infinitely many roots, and one usually disposes only of finitely many parameters that can be chosen in the controller design. If, however, the MID property holds and these parameters are chosen in such a way as to guarantee the existence of a root of maximal multiplicity, then this root is dominant, and hence determines the asymptotic behavior of the system, allowing for stabilization if one chooses this root with negative real part.

The aim of this paper is to start the investigation of whether the MID property holds for delay differential-algebraic systems of the form \eqref{IntroDiffAlgeb}. For this purpose, we restrict our attention to the first non-trivial situation, corresponding to the case $d_x = d_y = N = 1$ in which both equations in \eqref{IntroDiffAlgeb} are scalar and the system contains a single delay. This simple-looking, low-dimensional case illustrates many of the subtleties in the analysis of the MID property for delay differential-algebraic systems. Our main result for this system, Theorem~\ref{TheoDiffAlgeb}, shows that the MID property does hold in this setting, providing necessary and sufficient conditions on the system parameters for having a real root of maximal multiplicity and characterizing further the other roots of the characteristic quasipolynomial when these conditions are satisfied.

The paper is organized as follows. Section~\ref{SecMotivation} presents some examples of systems which can be put under the form \eqref{IntroDiffAlgeb}. We then present, in Section~\ref{SecRetarded}, a previous result from \cite{MazantiQualitative} on the MID property for retarded delay differential equations, which can be seen as a particular case of \eqref{IntroDiffAlgeb} in which $N = 1$ and $D_1 = 0$. We briefly recall the strategy of its proof, which serves as inspiration for the analysis of the MID property for a delay differential-algebraic system with scalar unknowns and a single delay in Section~\ref{SecDiffAlgeb}. An application to one of the examples from Section~\ref{SecMotivation} is provided in Section~\ref{SecAppli}.

\emph{Notation.} For a given complex number $s$, we denote by $\overline s$, $\Real s$, and $\Imag s$ its complex conjugate, real part, and imaginary part, respectively. Given nonnegative integers $n, k$ with $0 \leq k \leq n$, the notation $\binom{n}{k}$ represents the usual binomial coefficient $\frac{n!}{k! (n-k)!}$.

\section{Motivating examples}
\label{SecMotivation}

In this section, we briefly describe two systems that can be modeled by delay differential-algebraic equations.

\subsection{Electrical circuit with a transmission line}

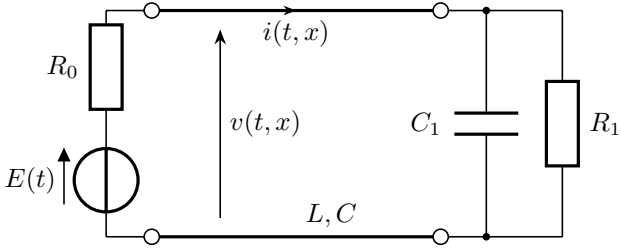
\begin{figure}[ht]
\centering
\begin{circuitikz}
\draw[semithick] (0.6, 0) circle[radius=0.1] (0.5, 0) -- (0, 0) to[V=$E(t)$] (0, 1.5) to[R=$R_0$] (0, 3) -- (0.5, 3) (0.6, 3) circle[radius=0.1];
\draw[very thick] (0.7, 3) -- (4.3, 3);
\draw[very thick] (0.7, 0) -- (4.3, 0);
\draw[semithick, -Stealth] (1.5, 0.25) -- node[midway, right] {$v(t, x)$} (1.5, 2.75);
\draw[-Stealth] (0.7, 3) -- (2.5, 3) node[below] {$i(t, x)$};
\draw[semithick] (4.4, 0) circle[radius=0.1] (4.5, 0) -- (5, 0) to[C=$C_1$] (5, 3) -- (4.5, 3) (4.4, 3) circle[radius=0.1];
\draw[semithick] (5, 3) -- (6, 3) to[R=$R_1$] (6, 0) -- (5, 0);
\fill (5, 0) circle[radius=0.05] (5, 3) circle[radius=0.05];
\draw (2.5, 0) node[above right] {$L, C$};
\end{circuitikz}
\caption{Electrical circuit with a transmission line}
\label{FigCircuit}
\end{figure}

Consider the electrical circuit from Figure~\ref{FigCircuit}, in which a voltage source $E(t)$ with internal resistance $R_0$ is connected, through a lossless transmission line of normalized length $1$ and normalized characteristic impedance $L$ and capacitance $C$, to a parallel association between a capacitor of capacitance $C_1$ and a resistor of resistance $R_1$. This electrical circuit can be described by the system
\begin{equation}
\label{ElectricCircuit}
\left\{
\begin{aligned}
& \partial_x v(t, x) + L \partial_t i(t, x) = 0, & & t \geq 0,\; x \in (0, 1), \\
& \partial_x i(t, x) + C \partial_t v(t, x) = 0, & & t \geq 0,\; x \in (0, 1), \\
& v(t, 0) = E(t) - R_0 i(t, 0), & & t \geq 0, \\
& i(t, 1) = C_1 \partial_t v(t, 1) + \frac{1}{R_1} v(t, 1), & & t \geq 0,
\end{aligned}
\right.
\end{equation}
with suitable initial conditions. Performing the classical change of variables into Riemann invariants
\[
\begin{aligned}
u_1(t, x) & = \frac{1}{2}\left[v(t, x) + \sqrt{\frac{L}{C}} i(t, x)\right], \\
u_2(t, x) & = \frac{1}{2}\left[v(t, x) - \sqrt{\frac{L}{C}} i(t, x)\right],
\end{aligned}
\]
and setting $y_1(t) = v(t, 1)$ and $y_2(t) = u_2(t, 1)$, one verifies that \eqref{ElectricCircuit} is rewritten as
\begin{equation}
\label{CircuitDiffAlgeb}
\left\{
\begin{aligned}
y_1^\prime(t) & = -\frac{1}{C_1} \left(\frac{1}{R_1} + \sqrt{\frac{C}{L}}\right) y_1(t) \\
 & - \frac{2}{C_1} \sqrt{\frac{C}{L}} \rho y_2(t - \tau) + \frac{2}{C_1} \sqrt{\frac{C}{L}} \frac{E(t - \tau/2)}{1 + R_0 \sqrt{\frac{L}{C}}}, \\
y_2(t) & = y_1(t) + \rho y_2(t - \tau) - \frac{E(t - \tau/2)}{1 + R_0 \sqrt{\frac{L}{C}}},
\end{aligned}
\right.
\end{equation}
where $\rho = \frac{1 - R_0 \sqrt{\frac{L}{C}}}{1 + R_0 \sqrt{\frac{L}{C}}}$ is the reflection coefficient at the voltage source and $\tau = 2 \sqrt{L C}$ represents the round-trip travel time along the transmission line. If one is interested in the stability analysis of \eqref{CircuitDiffAlgeb} under no input, i.e., when $E(t) = 0$ for every $t$, then \eqref{CircuitDiffAlgeb} reduces to a delay differential-algebraic system of the form \eqref{IntroDiffAlgeb} with $d_x = d_y = 1$ and $N = 1$. We refer the interested reader to \cite{Brayton1968Small, Halanay1997Stability, Niculescu2001Delay} and references therein for more details on these kind of circuits and their applications, as well as for further examples of engineering systems which can be put under the form of a delay differential-algebraic system \eqref{IntroDiffAlgeb}.

\subsection{Delayed output feedback for proper control systems}

Consider the linear control system
\begin{equation}
\label{LinControlSystem}
\left\{
\begin{aligned}
\dot x(t) & = A x(t) + B u(t), \\
y(t) & = C x(t) + D u(t),
\end{aligned}
\right.
\end{equation}
where $x(t) \in \mathbb R^{d_x}$ is the state, $u(t) \in \mathbb R^{d_u}$ is the control, $y(t) \in \mathbb R^{d_y}$ is the output, $d_x$, $d_u$, and $d_y$ are positive integers, and $A$, $B$, $C$, and $D$ are matrices of appropriate dimensions. We consider the problem of stabilizing \eqref{LinControlSystem} by a delayed output feedback of the form $u(t) = K y(t - \tau)$, where $\tau > 0$ is the delay and $K$ is a matrix of appropriate dimension. The closed-loop system is then
\begin{equation}
\label{DelayedOutputFeedback}
\left\{
\begin{aligned}
\dot x(t) & = A x(t) + B K y(t - \tau), \\
y(t) & = C x(t) + D K y(t - \tau),
\end{aligned}
\right.
\end{equation}
which is under the form \eqref{IntroDiffAlgeb} with $N = 1$ delay.

\section{Preliminary results on retarded delay differential equations}
\label{SecRetarded}

Consider the delay differential-algebraic system \eqref{IntroDiffAlgeb} in the particular case $N = 1$ and $D_1 = 0$. Letting $\tau = \tau_1$ and $B = B_1 C$, this system can be written under the form
\begin{equation}
\label{RetardedVector}
x^\prime(t) = A x(t) + B x(t - \tau).
\end{equation}
System \eqref{RetardedVector} is a system of delay differential equations said to be of \emph{retarded type} since the derivative of highest order only appears in the non-delayed term $x^\prime(t)$.

An important particular case which can be put into the form \eqref{RetardedVector} is that of a scalar retarded delay differential equation of order $n$,
\begin{equation}
\label{Retarded}
x^{(n)}(t) + \sum_{k=0}^{n-1} a_k x^{(k)}(t) + \sum_{k=0}^{n-1} \alpha_k x^{(k)}(t - \tau) = 0,
\end{equation}
where $n$ is a positive integer and the coefficients $a_k$ and $\alpha_k$ are real numbers for $k \in \{0, \dotsc, n - 1\}$. The corresponding characteristic quasipolynomial is the function $\Delta: \mathbb C \to \mathbb C$ given by
\begin{equation}
\label{RetardedDeltaN}
\Delta(s) = s^n + \sum_{k=0}^{n-1} a_k s^k + e^{-s\tau} \sum_{k=0}^{n-1} \alpha_k s^k.
\end{equation}
Notice that, according to Definition~\ref{DefiQuasipolynomial}, $\Delta$ is of degree $2n$, and hence any root of $\Delta$ has multiplicity at most $2n$.

The MID property has been extensively studied for some retarded differential equations under the form \eqref{Retarded}, such as in \cite{Boussaada2018Further, Boussaada2020Multiplicity, Boussaada2016Tracking, MazantiQualitative, MazantiSpectral}. In particular, \cite{MazantiQualitative} presents the following result.

\begin{thm}
\label{TheoRetarded}
Consider the quasipolynomial $\Delta$ from \eqref{RetardedDeltaN} and let $s_0 \in \mathbb R$.
\begin{enumerate}
\item\label{RetardedItemConditions} The number $s_0$ is a root of multiplicity $2n$ of $\Delta$ if and only if, for every $k \in \{0, \dotsc, n-1\}$,
\begin{equation}
\label{RetardedConditions}
\left\{
\begin{aligned}
a_k & = \binom{n}{k} (-s_0)^{n-k} \\
& \hphantom{=} {} + (-1)^{n-k} n! \sum_{j=k}^{n-1} \binom{j}{k} \binom{2n-j-1}{n-1} \frac{s_0^{j-k}}{j!\tau^{n-j}}, \\
\alpha_k & = (-1)^{n-1} e^{s_0 \tau} \sum_{j=k}^{n-1} \frac{(-1)^{j-k} (2n-j-1)!}{k! (j-k)! (n-j-1)!} \frac{s_0^{j-k}}{\tau^{n-j}}.
\end{aligned}
\right.
\end{equation}
\item\label{RetardedItemDominance} If \eqref{RetardedConditions} is satisfied, then $s_0$ is a strictly dominant root of $\Delta$.
\end{enumerate}
\end{thm}

The proof of this result is detailed in the case $n = 2$ in \cite{MazantiQualitative}, the full proof being provided in the extended version of that reference. A first step of the proof is to remark that it suffices to consider the case $s_0 = 0$ and $\tau = 1$, since the general case can be reduced to this setting by performing the translation and scaling of the spectrum represented by the change of variables $z = \tau(s - s_0)$. Part \ref{RetardedItemConditions} of Theorem~\ref{TheoRetarded} can then be obtained by straightforward computations, imposing that $\Delta^{(k)}(0) = 0$ for every $k \in \{0, \dotsc, 2n-1\}$. The dominance proof for establishing \ref{RetardedItemDominance} is then carried by providing first a priori bounds on the imaginary part of non-roots with non-negative real part and then using a suitable factorization of $\Delta$ to show, using this bound, that such roots cannot exist.

\section{MID for a delay differential-algebraic equation}
\label{SecDiffAlgeb}

We consider in this section the delay differential-algebraic equation
\begin{equation}
\label{DiffAlgeb}
\left\{
\begin{aligned}
x^\prime(t) & = a x(t) + b y(t - \tau), \\
y(t) & = c x(t) + d y(t - \tau),
\end{aligned}
\right.
\end{equation}
where $x(t) \in \mathbb R$, $y(t) \in \mathbb R$, and $a, b, c, d$ are real coefficients. System \eqref{DiffAlgeb} corresponds to \eqref{IntroDiffAlgeb} with $d_x = d_y = N = 1$ and the explicit computation of its characteristic quasipolynomial $\Delta$ from \eqref{IntroDelta} yields
\begin{equation}
\label{DiffAlgebDelta}
\Delta(s) = s - a - e^{-s\tau} \left(s d - ad + bc\right).
\end{equation}
Note that $\Delta$ is a quasipolynomial of degree $3$. The main result we prove in this paper is the following counterpart of Theorem~\ref{TheoRetarded}.

\begin{thm}
\label{TheoDiffAlgeb}
Consider the quasipolynomial $\Delta$ from \eqref{DiffAlgebDelta} and let $s_0 \in \mathbb R$.
\begin{enumerate}
\item\label{DiffAlgebItemConditions} The number $s_0$ is a root of multiplicity $3$ of $\Delta$ if and only if the coefficients $a, b, c, d$, the root $s_0$, and the delay $\tau$ satisfy the relations
\begin{equation}
\label{DiffAlgebConditions}
a = s_0 + \frac{2}{\tau}, \quad d = -e^{s_0 \tau}, \quad b c = -\frac{4}{\tau} e^{s_0 \tau}.
\end{equation}
\item\label{DiffAlgebItemDominance} If \eqref{DiffAlgebConditions} is satisfied, then $s_0$ is a dominant root of $\Delta$. Moreover, for every other complex root $s$ of $\Delta$, one has $\Real s = s_0$.
\item\label{DiffAlgebExplicit} Let $\Xi = \{\xi \in \mathbb R \mid \tan \xi = \xi\}$. If \eqref{DiffAlgebConditions} is satisfied, then the set of roots of $\Delta$ is $\{s_0 + i \frac{2}{\tau} \xi \suchthat \xi \in \Xi\}$.
\end{enumerate}
\end{thm}

\begin{rem}
With respect to Theorem~\ref{TheoRetarded}, Theorem~\ref{TheoDiffAlgeb} provides, in its part \ref{DiffAlgebExplicit}, additional information on the location of the other roots of $\Delta$. On the other hand, $s_0$ is not strictly dominant in this case.
\end{rem}

\begin{rem}
In the particular case $s_0 = 0$ and $\tau = 1$, \eqref{DiffAlgebConditions} yield $a = 2$, $d = -1$, and $b c = -4$. The corresponding quasipolynomial \eqref{DiffAlgebDelta} is then given by
\begin{equation}
\label{DiffAlgebDeltaHat}
\widehat\Delta(z) = z - 2 + e^{-z}(z + 2).
\end{equation}
For general $s_0 \in \mathbb R$ and $\tau > 0$, one may reduce to the above setting by performing the translation and scaling of the spectrum represented by the change of variables $z = \tau(s - s_0)$.
\end{rem}

The proof of Theorem~\ref{TheoDiffAlgeb} follows the same general line of that of Theorem~\ref{TheoRetarded}, but extra properties should be proved in order to obtain the additional conclusions of Theorem~\ref{TheoDiffAlgeb}. The main properties we need for the proof are provided in Appendix~\ref{SecAppendix}.

\begin{pf*}{Proof of Theorem~\ref{TheoDiffAlgeb}.}
Let $\widetilde\Delta$ be the quasipolynomial obtained from $\Delta$ by setting
\begin{equation}
\label{DefiTildeDelta}
\widetilde\Delta(z) = \tau\Delta\left(\frac{z}{\tau} + s_0\right)
\end{equation}
for $z \in \mathbb C$. Then
\[
\widetilde\Delta(z) = z + b_0 + e^{-z}(\beta_1 z + \beta_0)
\]
with
\begin{equation}
\label{RelationBA}
\begin{gathered}
b_0 = \tau (s_0 - a), \quad \beta_1 = - d e^{-s_0 \tau}, \\
\beta_0 =  \tau e^{-s_0 \tau} (a d - b c - d s_0).
\end{gathered}
\end{equation}
It follows immediately from relation \eqref{DefiTildeDelta} that $s_0$ is a root of multiplicity $3$ of $\Delta$ if and only if $0$ is a root of multiplicity $3$ of $\widetilde\Delta$. Since $\widetilde\Delta$ is a quasipolynomial of degree $3$, $0$ is a root of multiplicity $3$ of $\widetilde\Delta$ if and only if $\widetilde\Delta(0) = \widetilde\Delta^\prime(0) = \widetilde\Delta^{\prime\prime}(0) = 0$. We compute
\begin{align*}
\widetilde\Delta^\prime(z) & = 1 + e^{-z}(- \beta_1 z - \beta_0 + \beta_1), \\
\widetilde\Delta^{\prime\prime}(z) & = e^{-z}(\beta_1 z + \beta_0 - 2 \beta_1),
\end{align*}
and thus $0$ is a root of multiplicity $3$ of $\widetilde\Delta$ if and only if
\[
b_0 + \beta_0 = 0, \quad 1 - \beta_0 + \beta_1 = 0, \quad \beta_0 - 2 \beta_1 = 0.
\]
One immediately verifies that the above linear system of equations on $(b_0, \beta_1, \beta_0)$ admits a unique solution, given by $(b_0, \beta_1, \beta_0) = (-2, 1, 2)$. Using \eqref{RelationBA}, one concludes that $s_0$ is a root of multiplicity $3$ of $\Delta$ if and only if \eqref{DiffAlgebConditions} holds, concluding the proof of \ref{DiffAlgebItemConditions}. Notice moreover that, under \eqref{DiffAlgebConditions}, one has $\widetilde\Delta = \widehat\Delta$, where $\widehat\Delta$ is the quasipolynomial defined in \eqref{DiffAlgebDeltaHat}.

To prove \ref{DiffAlgebItemDominance}, it suffices to show that every root of $\widehat\Delta$ lies on the imaginary axis. Note first that
\begin{equation}
\label{Factorization}
\widehat\Delta(z) = z^3 \int_0^1 t (1-t) e^{-z t} \diff t,
\end{equation}
as one immediately verifies by integrating by parts. Assume, to obtain a contradiction, that there exists a root $z_0 \in \mathbb C$ of $\widehat\Delta$ such that $\Real z_0 \neq 0$. Writing $z_0 = \sigma + i \omega$ for $\sigma, \omega \in \mathbb R$ with $\sigma \neq 0$, one may assume, with no loss of generality thanks to Corollary~\ref{CoroSymmetry} in Appendix~\ref{SecAppendix}, that $\sigma > 0$ and $\omega > 0$. By Lemma~\ref{LemmAPrioriBound} in Appendix~\ref{SecAppendix}, one has $0 < \omega < 2$.

Using the fact that $z_0$ is a non-zero root of $\widehat\Delta$, one obtains from \eqref{Factorization} and taking the imaginary part that
\[
\int_0^1 t (1-t) e^{-\sigma t} \sin(\omega t) \diff t = 0.
\]
Since $0 < \omega < 2$, the function $t \mapsto t (1-t) e^{-\sigma t} \sin(\omega t)$ is strictly positive in $(0, 1)$, which contradicts the above equality. Hence \ref{DiffAlgebItemDominance} is proved.

Finally, \ref{DiffAlgebExplicit} follows immediately from the relation between $\widehat\Delta$ and $\Delta$ under \eqref{DiffAlgebConditions}, the fact that all roots of $\widehat\Delta$ lie on the imaginary axis, and Lemma~\ref{LemmImaginaryRoots} in Appendix~\ref{SecAppendix}.
\end{pf*}

\begin{rem}
With respect to other results on multiplicity-induced-dominancy for delay differential equations such as Theorem~\ref{TheoRetarded}, Theorem~\ref{TheoDiffAlgeb} provides, in its part \ref{DiffAlgebExplicit}, the additional information of the location of \emph{all} roots of $\Delta$. The set $\Xi$ of the real roots of the equation $\tan \xi = \xi$ is infinite, discrete, and can be written as $\Xi = \{\xi_k \suchthat k \in \mathbb Z\}$, where $(\xi_k)_{k \in \mathbb Z}$ is the increasing sequence of the roots of $\tan\xi = \xi$ with the convention that $\xi_0 = 0$. In particular, for every $k \in \mathbb Z$, one has $\xi_k \in \left(-\frac{\pi}{2} + k \pi, \frac{\pi}{2} + k \pi\right)$ and $\xi_{-k} = - \xi_k$. We also recall that $\xi_k = k \pi + \frac{\pi}{2} + o(1)$ as $k \to +\infty$.
\end{rem}

\section{Application to the delayed output feedback for proper control systems}
\label{SecAppli}

Consider system \eqref{DelayedOutputFeedback} with $d_x = 1$, $d_y = 2$, $d_u = 1$, which we write as
\begin{equation}
\label{AppliSystem}
\left\{
\begin{aligned}
\dot x(t) & = a x(t) + b k_1 y_1(t - \tau) + b k_2 y_2(t - \tau), \\
y_1(t) & = c_1 x(t) + d_1 k_1 y_1(t - \tau) + d_1 k_2 y_2(t - \tau), \\
y_2(t) & = c_2 x(t) + d_2 k_1 y_1(t - \tau) + d_2 k_2 y_2(t - \tau). \\
\end{aligned}
\right.
\end{equation}
From \eqref{IntroDelta}, we compute its characteristic quasipolynomial
\begin{multline}
\label{AppliDelta}
\Delta(s) = s - a - e^{-s\tau}\bigl((d_1 k_1 + d_2 k_2)s \\ - k_1(a d_1 - b c_1) - k_2(a d_2 - b c_2)\bigr).
\end{multline}
Even though \eqref{AppliSystem} is not under the form \eqref{DiffAlgeb}, its characteristic quasipolynomial \eqref{AppliDelta} is of the same form of that of \eqref{DiffAlgebDelta}, and thus Theorem~\ref{TheoDiffAlgeb} can be applied to \eqref{AppliDelta}. We wish to design the parameters $k_1$ and $k_2$ and obtain conditions on the delay $\tau$ in order to achieve maximal multiplicity of some root $s_0 < 0$. Then Theorem~\ref{TheoDiffAlgeb} will ensure the dominance of this root, implying the exponential stability of the system.

Conditions \eqref{DiffAlgebConditions} can be rewritten in the context of \eqref{AppliDelta} as
\begin{equation}
\label{AppliConditions}
\begin{aligned}
& a = s_0 + \frac{2}{\tau}, \\
& d_1 k_1 + d_2 k_2 = - e^{s_0 \tau}, \\
& (a d_1 - b c_1) k_1 + (a d_2 - b c_2) k_2 = \left(\frac{2}{\tau} - s_0\right) e^{s_0 \tau}.
\end{aligned}
\end{equation}
Note that the first equation in \eqref{AppliConditions} can be rewritten as $s_0 = a - \frac{2}{\tau}$, hence one may stabilize the system by a root of maximal multiplicity only if $a \leq 0$ or $\tau < \frac{2}{a}$.

Let us denote $\delta_i = a d_i - b c_i$ for $i \in \{1, 2\}$. The second and third equations of \eqref{AppliConditions} can be rewritten as
\begin{equation}
\label{SystK}
\begin{pmatrix}d_1 & d_2 \\ \delta_1 & \delta_2 \\\end{pmatrix} \begin{pmatrix}k_1 \\ k_2 \\\end{pmatrix} = \begin{pmatrix}-1 \\ \tfrac{4}{\tau} - a\end{pmatrix} e^{a \tau - 2}.
\end{equation}
This system admits at least one solution if and only if $d_1 \delta_2 \neq d_2 \delta_1$ or $\delta_i = \left(a - \frac{4}{\tau}\right) d_i$ for $i \in \{1, 2\}$, with exactly one solution in the first case and infinitely many solutions in the second case. This discussion can be concentrated in the following result.

\begin{prop}
\label{PropAppli}
Consider system \eqref{AppliSystem} for given real parameters $a$, $b$, $c_1$, $c_2$, $d_1$, and $d_2$ and a positive delay $\tau$, and assume that either $a \leq 0$ or $\tau < \frac{2}{a}$. Let $\delta_i = a d_i - b c_i$ for $i \in \{1, 2\}$ and assume moreover that either $d_1 \delta_2 \neq d_2 \delta_1$ or $\delta_i = \left(a - \frac{4}{\tau}\right) d_i$ for $i \in \{1, 2\}$. Then there exist real parameters $k_1$, $k_2$ such that \eqref{AppliSystem} is exponentially stable, with exponential decay rate $a - \frac{2}{\tau}$. Moreover, $k_1$ and $k_2$ are solutions of the linear system \eqref{SystK}.
\end{prop}

\begin{figure}[ht]
\centering
\resizebox{\columnwidth}{!}{\input{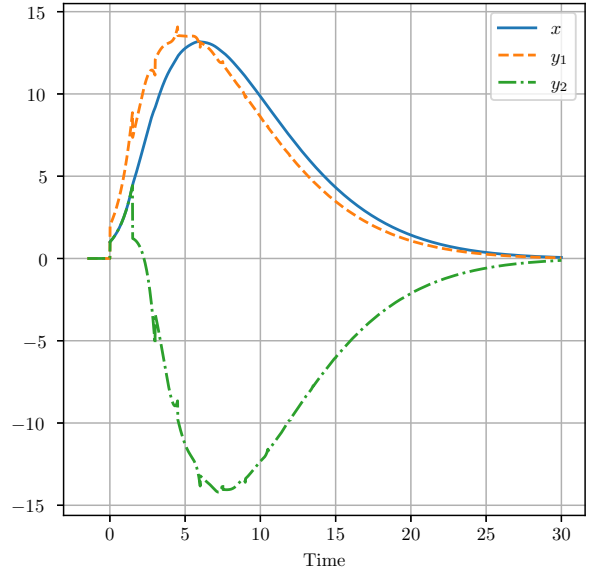}}
\caption{Trajectories of \eqref{AppliSystem} with the designed delayed output feedback.}
\label{AppliFigure}
\end{figure}

Consider, as an example, the case $a = 1$, $b = 1$, $c_1 = 2$, $c_2 = 1$, $d_1 = 1$, $d_2 = 2$, and $\tau = \frac{3}{2}$. Since $a = 1$, the corresponding open-loop system is unstable. Note that the inequality $\tau < \frac{2}{a}$ is indeed satisfied and the exponential decay rate one may obtain with Proposition~\ref{PropAppli} is $a - \frac{2}{\tau} = -\frac{1}{3}$. One computes, using Proposition~\ref{PropAppli}, the feedback parameters $k_1 \approx -0.87610$ and $k_2 \approx 0.13478$. Figure~\ref{AppliFigure} presents a numerical simulation of the solutions of \eqref{AppliConditions} with these parameters and with initial conditions $x(0) = 1$ and $y_1(t) = y_2(t) = 0$ for $t < 0$. One observes, as expected, that solutions converge to $0$ exponentially.

\begin{ack}
This work is supported by a public grant overseen by the French National Research Agency (ANR) as part of the ``Investissement d'Avenir'' program, through the iCODE project funded by the IDEX Paris-Saclay, ANR-11-IDEX-0003-02. The authors also acknowledge the support of Institut Polytechnique des Sciences Avanc\'ees (IPSA).
\end{ack}

\bibliography{Bib}

\begin{thebibliography}{28}
\providecommand{\natexlab}[1]{#1}
\providecommand{\url}[1]{\texttt{#1}}
\providecommand{\urlprefix}{URL }
\expandafter\ifx\csname urlstyle\endcsname\relax
  \providecommand{\doi}[1]{doi:\discretionary{}{}{}#1}\else
  \providecommand{\doi}{doi:\discretionary{}{}{}\begingroup
  \urlstyle{rm}\Url}\fi

\bibitem[{Bellman and Cooke(1963)}]{Bellman1963Differential}
Bellman, R. and Cooke, K.L. (1963).
\newblock \emph{Differential-difference equations}.
\newblock Academic Press, New York-London.

\bibitem[{Boussaada and
  Niculescu(2016{\natexlab{a}})}]{Boussaada2016Characterizing}
Boussaada, I. and Niculescu, S.I. (2016{\natexlab{a}}).
\newblock Characterizing the codimension of zero singularities for time-delay
  systems: a link with {V}andermonde and {B}irkhoff incidence matrices.
\newblock \emph{Acta Appl. Math.}, 145, 47--88.

\bibitem[{Boussaada and Niculescu(2016{\natexlab{b}})}]{Boussaada2016Tracking}
Boussaada, I. and Niculescu, S.I. (2016{\natexlab{b}}).
\newblock Tracking the algebraic multiplicity of crossing imaginary roots for
  generic quasipolynomials: a {V}andermonde-based approach.
\newblock \emph{IEEE Trans. Automat. Control}, 61(6), 1601--1606.

\bibitem[{Boussaada et~al.(2020)Boussaada, Niculescu, El~Ati, Perez-Ramos, and
  Trabelsi}]{Boussaada2020Multiplicity}
Boussaada, I., Niculescu, S.I., El~Ati, A., Perez-Ramos, R., and Trabelsi, K.L.
  (2020).
\newblock Multiplicity-induced-dominancy in parametric second-order delay
  differential equations: Analysis and application in control design.
\newblock \emph{ESAIM Control Optim. Calc. Var.}

\bibitem[{Boussaada et~al.(2018)Boussaada, Tliba, Niculescu, \"{U}nal, and
  Vyhl\'{\i}dal}]{Boussaada2018Further}
Boussaada, I., Tliba, S., Niculescu, S.I., \"{U}nal, H.U., and Vyhl\'{\i}dal,
  T. (2018).
\newblock Further remarks on the effect of multiple spectral values on the
  dynamics of time-delay systems. {A}pplication to the control of a mechanical
  system.
\newblock \emph{Linear Algebra Appl.}, 542, 589--604.

\bibitem[{Brayton(1968)}]{Brayton1968Small}
Brayton, R.K. (1968).
\newblock Small-signal stability criterion for electrical networks containing
  lossless transmission lines.
\newblock \emph{{IBM} Journal of Research and Development}, 12(6), 431--440.

\bibitem[{Brenan et~al.(1989)Brenan, Campbell, and
  Petzold}]{Brenan1989Numerical}
Brenan, K.E., Campbell, S.L., and Petzold, L.R. (1989).
\newblock \emph{Numerical solution of initial value problems in
  differential-algebraic equations}.
\newblock North-Holland Publishing Co., New York.

\bibitem[{Coleman(1998)}]{Coleman1998Differential}
Coleman, R. (ed.) (1998).
\newblock \emph{Differential algebraic equations}.
\newblock Baltzer Science Publishers BV, Bussum.
\newblock Selected papers from the International Symposium held in Grenoble,
  May 1997, Numer. Algorithms {{\bf{1}}9} (1998), no. 1-4.

\bibitem[{Cooke and van~den Driessche(1986)}]{Cooke1986Zeroes}
Cooke, K.L. and van~den Driessche, P. (1986).
\newblock On zeroes of some transcendental equations.
\newblock \emph{Funkcial. Ekvac.}, 29(1), 77--90.

\bibitem[{Diekmann et~al.(1995)Diekmann, van Gils, Verduyn~Lunel, and
  Walther}]{Diekmann1995Delay}
Diekmann, O., van Gils, S.A., Verduyn~Lunel, S.M., and Walther, H.O. (1995).
\newblock \emph{Delay equations: Functional-, complex-, and nonlinear
  analysis}, volume 110 of \emph{Applied Mathematical Sciences}.
\newblock Springer-Verlag, New York.

\bibitem[{Griepentrog and M\"{a}rz(1986)}]{Griepentrog1986Differential}
Griepentrog, E. and M\"{a}rz, R. (1986).
\newblock \emph{Differential-algebraic equations and their numerical
  treatment}, volume~88 of \emph{Teubner-Texte zur Mathematik [Teubner Texts in
  Mathematics]}.
\newblock BSB B. G. Teubner Verlagsgesellschaft, Leipzig.

\bibitem[{Gu et~al.(2003)Gu, Kharitonov, and Chen}]{Gu2003Stability}
Gu, K., Kharitonov, V.L., and Chen, J. (2003).
\newblock \emph{Stability of time-delay systems}.
\newblock Control Engineering. Birkh\"{a}user Boston, Inc., Boston, MA.

\bibitem[{Halanay(1966)}]{Halanay1966Differential}
Halanay, A. (1966).
\newblock \emph{Differential equations: {S}tability, oscillations, time lags}.
\newblock Academic Press, New York-London.

\bibitem[{Halanay and Rasvan(1997)}]{Halanay1997Stability}
Halanay, A. and Rasvan, V. (1997).
\newblock Stability radii for some propagation models.
\newblock \emph{IMA J. Math. Control Inform.}, 14(1), 95--107.
\newblock Distributed parameter systems: analysis, synthesis and applications,
  Part 1.

\bibitem[{Hale and Verduyn~Lunel(1993)}]{Hale1993Introduction}
Hale, J.K. and Verduyn~Lunel, S.M. (1993).
\newblock \emph{Introduction to functional-differential equations}, volume~99
  of \emph{Applied Mathematical Sciences}.
\newblock Springer-Verlag, New York.

\bibitem[{Insperger and St\'{e}p\'{a}n(2011)}]{Insperger2011Semi}
Insperger, T. and St\'{e}p\'{a}n, G. (2011).
\newblock \emph{Semi-discretization for time-delay systems}, volume 178 of
  \emph{Applied Mathematical Sciences}.
\newblock Springer, New York.
\newblock Stability and engineering applications.

\bibitem[{Kumar and Daoutidis(1999)}]{Kumar1999Control}
Kumar, A. and Daoutidis, P. (1999).
\newblock \emph{Control of nonlinear differential algebraic equation systems},
  volume 397 of \emph{Chapman \& Hall/CRC Research Notes in Mathematics}.
\newblock Chapman \& Hall/CRC, Boca Raton, FL.
\newblock With applications to chemical processes.

\bibitem[{Kunkel and Mehrmann(2006)}]{Kunkel2006Differential}
Kunkel, P. and Mehrmann, V. (2006).
\newblock \emph{Differential-algebraic equations}.
\newblock EMS Textbooks in Mathematics. European Mathematical Society (EMS),
  Z\"{u}rich.
\newblock Analysis and numerical solution.

\bibitem[{Mazanti et~al.(2020{\natexlab{a}})Mazanti, Boussaada, and
  Niculescu}]{MazantiQualitative}
Mazanti, G., Boussaada, I., and Niculescu, S.I. (2020{\natexlab{a}}).
\newblock On qualitative properties of single-delay linear retarded
  differential equations: Characteristic roots of maximal multiplicity are
  necessarily dominant.
\newblock Preprint.

\bibitem[{Mazanti et~al.(2020{\natexlab{b}})Mazanti, Boussaada, Niculescu, and
  Vyhl\'{\i}dal}]{MazantiSpectral}
Mazanti, G., Boussaada, I., Niculescu, S.I., and Vyhl\'{\i}dal, T.
  (2020{\natexlab{b}}).
\newblock Spectral dominance of complex roots for single-delay linear
  equations.
\newblock Preprint.

\bibitem[{Michiels and Niculescu(2007)}]{Michiels2007Stability}
Michiels, W. and Niculescu, S.I. (2007).
\newblock \emph{Stability and sta\-bi\-li\-za\-tion of time-delay systems: An
  ei\-gen\-va\-lue-based approach}, volume~12 of \emph{Advances in Design and
  Control}.
\newblock SIAM, Philadelphia, PA.

\bibitem[{Mori et~al.(1982)Mori, Fukuma, and Kuwahara}]{Mori1982Estimate}
Mori, T., Fukuma, N., and Kuwahara, M. (1982).
\newblock On an estimate of the decay rate for stable linear delay systems.
\newblock \emph{Internat. J. Control}, 36(1), 95--97.

\bibitem[{Niculescu(2001)}]{Niculescu2001Delay}
Niculescu, S.I. (2001).
\newblock \emph{Delay effects on stability: A robust control approach}, volume
  269 of \emph{Lecture Notes in Control and Information Sciences}.
\newblock Springer-Verlag London Ltd., London.

\bibitem[{Olgac and Sipahi(2002)}]{Olgac2002Exact}
Olgac, N. and Sipahi, R. (2002).
\newblock An exact method for the stability analysis of time-delayed linear
  time-invariant ({LTI}) systems.
\newblock \emph{IEEE Trans. Automat. Control}, 47(5), 793--797.

\bibitem[{Pinney(1958)}]{Pinney1958Ordinary}
Pinney, E. (1958).
\newblock \emph{Ordinary difference-differential equations}.
\newblock University of California Press, Berkeley-Los Angeles.

\bibitem[{P\'{o}lya and Szeg\H{o}(1998)}]{Polya1998Problems}
P\'{o}lya, G. and Szeg\H{o}, G. (1998).
\newblock \emph{Problems and theorems in analysis. {I}}.
\newblock Classics in Mathematics. Springer-Verlag, Berlin.
\newblock Series, integral calculus, theory of functions, Translated from the
  German by Dorothee Aeppli, Reprint of the 1978 English translation.

\bibitem[{Sipahi et~al.(2011)Sipahi, Niculescu, Abdallah, Michiels, and
  Gu}]{Sipahi2011Stability}
Sipahi, R., Niculescu, S.I., Abdallah, C.T., Michiels, W., and Gu, K. (2011).
\newblock Stability and stabilization of systems with time delay: limitations
  and opportunities.
\newblock \emph{IEEE Control Syst. Mag.}, 31(1), 38--65.

\bibitem[{St\'{e}p\'{a}n(1989)}]{Stepan1989Retarded}
St\'{e}p\'{a}n, G. (1989).
\newblock \emph{Retarded dynamical systems: stability and characteristic
  functions}, volume 210 of \emph{Pitman Research Notes in Mathematics Series}.
\newblock Longman Scientific \& Technical, Harlow; copublished in the United
  States with John Wiley \& Sons, Inc., New York.

\end{thebibliography}

\appendix

\section{Some technical results}
\label{SecAppendix}

We present in this appendix technical results used in the proof of Theorem~\ref{TheoDiffAlgeb}. We start by providing some properties of the quasipolynomial $\widehat\Delta$ from \eqref{DiffAlgebDeltaHat}. The first one is the following identity, whose proof is straightforward.

\begin{lem}
Let $\widehat\Delta$ be given by \eqref{DiffAlgebDeltaHat}. Then, for every $z \in \mathbb C$, one has
\[\widehat\Delta(-z) = -e^{z} \widehat\Delta(z).\]
\end{lem}

As a consequence of the previous identity, one immediately obtains the following symmetry property of the roots of $\widehat\Delta$.

\begin{cor}
\label{CoroSymmetry}
Let $\widehat\Delta$ be given by \eqref{DiffAlgebDeltaHat} and assume that $z_0 \in \mathbb C$ is such that $\widehat\Delta(z_0) = 0$. Then $\widehat\Delta(z_0) = \widehat\Delta(\overline z_0) = \widehat\Delta(-z_0) = \widehat\Delta(-\overline z_0) = 0$.
\end{cor}

We also need the following a priori bound on the imaginary part of the roots of $\widehat\Delta$ outside of the imaginary axis.

\begin{lem}
\label{LemmAPrioriBound}
Let $\widehat\Delta$ be given by \eqref{DiffAlgebDeltaHat} and assume that $z_0 \in \mathbb C$ is such that $\Real z_0 \neq 0$ and $\widehat\Delta(z_0) = 0$. Then $\abs{\Imag z_0} < 2$.
\end{lem}

\begin{pf}
Let $z_0 \in \mathbb C$ be as in the statement. We write $z_0 = \sigma + i \omega$ with $\sigma, \omega \in \mathbb R$ and $\sigma \neq 0$. Thanks to Corollary~\ref{CoroSymmetry}, we assume, with no loss of generality, that $\sigma > 0$. Since $z_0$ is a root of $\widehat\Delta$, one has $e^{- z_0} (z_0 + 2) = 2 - z_0$ and thus, in particular, $\abs{z_0 + 2}^2 = e^{2 \sigma} \abs{2 - z_0}^2$, which yields $(\sigma + 2)^2 + \omega^2 = e^{2 \sigma} \left((2 - \sigma)^2 + \omega^2\right)$.

Let $f: \mathbb R \to \mathbb R$ be the function defined by $f(x) = e^{2 x} \left((2 - x)^2 + \omega^2\right) - (x + 2)^2 - \omega^2$. Note that $f(0) = f(\sigma) = 0$. Since $f$ is differentiable and $\sigma > 0$, the mean value theorem yields the existence of $x_\ast \in (0, \sigma)$ such that $f^\prime(x_\ast) = 0$. We compute $f^\prime(x) = 2 e^{2x} \left((2-x)(1-x) + \omega^2\right) - 2(x+2)$, and thus
\[
e^{2x_\ast} \left((2-x_\ast)(1-x_\ast) + \omega^2\right) = x_\ast + 2.
\]
Since one has further that $e^{2x_\ast} > 1$, one deduces that
\[
(2-x_\ast)(1-x_\ast) + \omega^2 < x_\ast + 2,
\]
which is equivalent to
\[
x_\ast^2 - 4 x_\ast + \omega^2 < 0.
\]
Letting $g: \mathbb R \to \mathbb R$ be the polynomial $g(x) = x^2 - 4 x + \omega^2$, since $\lim_{x \to \pm\infty} g(x) = +\infty$, the above inequality implies that $g$ must admit two distinct real roots, and thus its discriminant is positive, i.e., $16 - 4 \omega^2 > 0$, which is equivalent to $\omega^2 < 4$. Thus $\abs{\omega} < 2$, as required.
\end{pf}

As a final technical result, we provide the following characterization of the roots of $\widehat\Delta$ on the imaginary axis.

\begin{lem}
\label{LemmImaginaryRoots}
Let $\widehat\Delta$ be given by \eqref{DiffAlgebDeltaHat}, $\Xi$ be as in the statement of Theorem~\ref{TheoDiffAlgeb}\ref{DiffAlgebExplicit}, and $\zeta \in \mathbb R$. Then $i\zeta$ is a root of $\widehat\Delta$ if and only if $\frac{\zeta}{2} \in \Xi$.
\end{lem}

\begin{pf}
Note that $i\zeta$ is a root of $\widehat\Delta$ if and only if
\[
i\zeta - 2 + e^{-i\zeta} (i\zeta + 2) = 0,
\]
which is the case if and only if
\[
\left\{
\begin{aligned}
-2 + 2 \cos\zeta + \zeta\sin\zeta & = 0, \\
\zeta + \zeta \cos\zeta - 2 \sin\zeta & = 0.
\end{aligned}
\right.
\]
The above system is equivalent to
\begin{equation}
\label{ExplicitSystem}
R_{-\zeta} \begin{pmatrix}2 \\ \zeta\end{pmatrix} = \begin{pmatrix}2 \\ -\zeta\end{pmatrix},
\end{equation}
where, for $\theta \in \mathbb R$, $R_\theta$ is the rotation matrix in $\mathbb R^2$, defined by
\[
R_\theta = \begin{pmatrix}\cos\theta & -\sin\theta \\ \sin\theta & \cos\theta\end{pmatrix}.
\]
Recalling that $R_{\theta}^{-1} = R_{-\theta}$ and $R_{\theta_1} R_{\theta_2} = R_{\theta_1 + \theta_2}$ for every $(\theta, \theta_1, \theta_2) \in \mathbb R^3$, one deduces that \eqref{ExplicitSystem} is equivalent to
\[
R_{-\frac{\zeta}{2}} \begin{pmatrix}1 \\ \frac{\zeta}{2}\end{pmatrix} = R_{\frac{\zeta}{2}} \begin{pmatrix}1 \\ -\frac{\zeta}{2}\end{pmatrix}.
\]
One then immediately verifies that the above system is equivalent to
\[
-\sin\left(\frac{\zeta}{2}\right) + \frac{\zeta}{2} \cos\left(\frac{\zeta}{2}\right) = \sin\left(\frac{\zeta}{2}\right) - \frac{\zeta}{2} \cos\left(\frac{\zeta}{2}\right),
\]
which holds if and only if
\[
\tan\left(\frac{\zeta}{2}\right) = \frac{\zeta}{2},
\]
i.e., if and only if $\frac{\zeta}{2} \in \Xi$, as required.
\end{pf}

\end{document}